\documentclass{article}%
\usepackage{amsmath}
\usepackage{amsfonts}
\usepackage{amssymb}
\usepackage{graphicx}%
\providecommand{\U}[1]{\protect\rule{.1in}{.1in}}

\ifx\pdfoutput\relax\let\pdfoutput=\undefined\fi
\newcount\msipdfoutput
\ifx\pdfoutput\undefined\else
\ifcase\pdfoutput\else
\ifx\paperwidth\undefined\else
\ifdim\paperheight=0pt\relax\else\pdfpageheight\paperheight\fi
\ifdim\paperwidth=0pt\relax\else\pdfpagewidth\paperwidth\fi
\fi\fi\fi
\begin{document}

\title{How to Measure Statistical Evidence and Its Strength: Bayes Factors or
Relative Belief Ratios?}
\author{L. Al-Labadi, A. Alzaatreh and M. Evans\\U. of Toronto, American U. of Sharjah and U. of Toronto}
\date{}
\maketitle

\begin{abstract}
Both the Bayes factor and the relative belief ratio satisfy the principle of
evidence and so can be seen to be valid measures of statistical evidence.
Certainly Bayes factors are regularly employed. The question then is: which of
these measures of evidence is more appropriate? It is argued here that there
are questions concerning the validity of a current commonly used definition of
the Bayes factor based on a mixture prior and, when all is considered, the
relative belief ratio has better properties as a measure of evidence. It is
further shown that, when a natural restriction on the mixture prior is
imposed, the Bayes factor equals the relative belief ratio obtained without
using the mixture prior. Even with this restriction, this still leaves open
the question of how the strength of evidence is to be measured. It is argued
here that the current practice of using the size of the Bayes factor to
measure strength is not correct and a solution to this issue is presented.
Several general criticisms of these measures of evidence are also discussed
and addressed.

\end{abstract}

\bigskip


\section{Introduction}

One of the virtues of the Bayesian approach to statistical analysis is that it
gives an unambiguous definition of what it means for there to be evidence in
favor of or against a particular value of a parameter. This is provided by the
following principle.

\begin{quote}
\textit{Principle of Evidence}: if the posterior probability of an event is
greater than (less than, equal to) its prior probability, then there is
evidence in favor of (against, no evidence either way of) the event being true.
\end{quote}

\noindent This seems like a very simple and intuitively satisfying way of
characterizing evidence and it has long been considered to be quite natural
and obvious. For example,

\begin{quote}
Popper (1968) The Logic of Scientific Discovery, Appendix ix \textquotedblleft
If we are asked to give a criterion of the fact that the evidence $y$ supports
or corroborates a statement $x$, the most obvious reply is: that $y$ increases
the probability of $x.$"

Achinstein (2001) \textquotedblleft for a fact $e$ to be evidence that a
hypothesis $h$ is true, it is both necessary and sufficient for $e$ to
increase $h$'s probability over its prior probability\textquotedblright.
\end{quote}

\noindent Both of these references are by philosophers of science and indeed
there is much discussion of this idea and its consequences in that literature
where it is often referred to as confirmation theory, see Salmon (1973) and
the Stanford Encyclopedia of Philosophy (2020).

While the statistical literature has made attempts to deal directly with the
concept of statistical evidence, Birnbaum (1964) and Royall (1997) being
notable examples, it is fair to say that the principle of evidence has not had
a large influence on the development of statistical methodology. This is
unfortunate because the ingredients that are available in a statistical
analysis can result in many of the issues raised about confirmation theory by
philosophers being successfully addressed and new and interesting inferences
with many good properties arise. Evans (2015) contains a development of a
statistical theory of inference based largely on this principle.

The \textit{Bayes factor} is however a notable exception. For this suppose we
have a sampling model $(\mathcal{X},\mathcal{B},\{P_{\theta}:\theta\in
\Theta\})$ and a prior probability model $(\Theta,\mathcal{A},\Pi).$ Suppose
that $A\in\mathcal{A}$ is an event whose truth value it is desired to know and
suppose further that $0<\Pi(A)<1.$ Let $\Pi(A\,|\,x)$ denote the posterior
probability of $A$ after observing the data $x\in\mathcal{X}.$ The
\textit{Bayes factor} in favor of $A$ being true is then given by
\[
BF(A\,|\,x)=\left\{  \frac{\Pi(A\,|\,x)}{\Pi(A^{c}\,|\,x)}\right\}  \left\{
\frac{\Pi(A)}{\Pi(A^{c})}\right\}  ^{-1},
\]
namely, the ratio of the posterior odds in favor of $A$ to the prior odds in
favor of $A.$ The following result is then an immediate application of the
principle of evidence.\smallskip

\noindent\textbf{Proposition 1. }If $0<\Pi(A)<1,$ then there is evidence in
favor of (against, no evidence either way of) $A$ being true iff
$BF(A\,|\,x)>(<,=)1.$\smallskip

Any measure of evidence that satisfies the principle of evidence using some
cut-off, as the Bayes factor does with cut-off 1, will be called a
\textit{valid} measure of evidence. Note that the meaning of the equality of
any measure of evidence to its cut-off simply means that the events
$A\subset\Theta$ and $\{x\}\subset\mathcal{X}$ are statistically independent
and so, at least in the discrete case, observing $x$ can tell us nothing about
the truth of $A.$ Also, any 1-1 function of a valid measure of evidence can be
used instead, such as $\log BF(A\,|\,x)$ with cut-off $0,$ as this is also a
valid measure of evidence, see Section 2.4.

One immediate question about the Bayes factor is why is it necessary to use
the ratio of odds? For example, a somewhat simpler measure of evidence is
given by the \textit{relative belief ratio}
\[
RB(A\,|\,x)=\frac{\Pi(A\,|\,x)}{\Pi(A)},
\]
namely, the ratio of the posterior to prior probability of $A,$ as this is
also a valid measure of evidence again with the cut-off 1. Of course, there
are still other valid measures of evidence, such as the difference
$\Pi(A\,|\,x)-\Pi(A)$ with cut-off 0, but there are some natural reasons to
focus on the Bayes factor and relative belief ratio that will be discussed subsequently.

When $\Pi$ is discrete with density $\pi$ with respect to counting measure and
$\Pi(A)=p,$ then%
\[
BF(A\,|\,x)=\frac{%
{\displaystyle\sum\limits_{\theta\in A}}
\pi(\theta)f_{\theta}(x)}{%
{\displaystyle\sum\limits_{\theta\in A^{c}}}
\pi(\theta)f_{\theta}(x)}\frac{1-p}{p}=\frac{m(x\,|\,A)}{m(x\,|\,A^{c})}%
\]
where $m(\cdot\,|\,A),m(\cdot\,|\,A^{c})$ denote the prior predictive
densities for $x$ given $\theta\in A$ and $\theta\in A^{c},$ respectively. By
contrast
\begin{equation}
RB(A\,|\,x)=\frac{m(x\,|\,A)}{m(x)}, \label{reb1}%
\end{equation}
where $m$ is the full prior predictive. These formulas can lead to considering
the Bayes factor as something like a likelihood ratio and the relative belief
ratio as something like a normalized likelihood but these interpretations are
definitely subsidiary to their roles as valid measures of evidence, however,
see Section 2.4.

The use of Bayes factors in inference, and in particular for hypothesis
assessment, is quite common but the relative belief ratio is not used to
anywhere near the same extent. In fact, it seems that some believe that\ there
are natural advantages to the Bayes factor. It is the purpose of this paper to
show that the true state of affairs is exactly the opposite. In fact, we
believe that using relative belief ratios as the basis for Bayesian inference
simplifies the development of inference based on a measure of evidence dramatically.

One advantage that is sometimes claimed for the Bayes factor, is that its
value also measures the strength of the evidence. So, for example, the bigger
the value of $BF(A\,|\,x)>1$ is, the stronger the evidence is in favor of $A$
being true and, of course, the same could be claimed for $RB(A\,|\,x).$ While
there is an element of truth in such claims, a substantial issue needs to be
addressed and this is typically ignored, namely, is there a universal scale on
which evidence can be measured? This is discussed in Section 3.2 and it will
be argued that the appropriate approach to assessing the strength of evidence
is to measure \textit{how strongly it is believed what the evidence says}. In
other words, strength of evidence must be measured by a posterior probability.

One may wonder why it is necessary to base inferences on a measure of evidence
as opposed to the commonly used framework where a loss function is employed.
This issue is not fully discussed here but note the appeal of not having to
choose yet another ingredient, the loss function, to obtain inferences.
Perhaps even more important is the fact that the ingredients to an
evidence-based approach, namely, the model and prior, can both be checked
against the data to see if there is a problem with these ingredients, via
model checking and checking for prior-data conflict, respectively, while such
a check does not seem to be generally available for loss functions. The
importance\ of this checking arises from the scientific principle that all
subjective choices to an analysis should be falsifiable and indeed, need to be
checked against the objective, at least when it is collected correctly, data
and replaced when found wanting. Without such checking any
inferences/conclusions drawn based on the ingredients are suspect as perhaps
these are chosen to achieve some specific goal.\ Some relevant references for
checking for prior-data conflict are Evans and Moshonov (2006), Evans and Jang
(2011a), Nott et al. (2021) while Evans and Jang (2011b) discusses what can be
done when such a conflict arises.

In Section 2 a comparison is made between Bayes factors and relative belief
ratios. In Section 3 some criticisms of these measures are addressed and, in
particular, the issue of assessing the strength of evidence. Section 4 draws
some overall conclusions. While parts of the paper comprise a review of known
properties of Bayes factors and relative belief ratios, there is also a great
deal that is new. For example, Propositions 5-10 and 12 are all new results.

\section{Bayes Factor versus Relative Belief Ratio}

A comparison between the relative belief ratio and the Bayes factor is now
considered in the following subsections.

\subsection{Simplicity}

Certainly it seems obvious that the relative belief ratio is a simpler
function of the ingredients than the Bayes factor. In particular, the Bayes
factor is a simple function of relative belief ratios alone but not
conversely.\smallskip\ 

\noindent\textbf{Proposition 2. }If $0<\Pi(A)<1,$ then
$BF(A\,|\,x)=RB(A\,|\,x)/RB(A^{c}\,|\,x).$\smallskip

\noindent From this we see that the Bayes factor is a comparison of the
evidence for $A$ with the evidence for $A^{c}.$ But this comparison doesn't
really make sense because it is immediate that $RB(A\,|\,x)>(<,=)1$ iff
$RB(A^{c}\,|\,x)<(>,=)1.$ In other words, if by the principle of evidence
there is evidence in favor $A,$ then there is always evidence against $A^{c}$
and conversely. It is also immediate from Proposition 2 that, if there is
evidence in favor of $A,$ then $BF(A\,|\,x)\geq RB(A\,|\,x)$ and if there is
evidence against $A,$ then $BF(A\,|\,x)\leq RB(A\,|\,x).$ As discussed in
Section 3.2, however, this cannot be interpreted as the Bayes factor
indicating stronger evidence, either in favor or against, when compared to the
relative belief ratio. The main point in this section is that there is no
particular reason to consider the ratio of odds in favor instead of the ratio
of the probabilities.

\subsection{Events with Prior Content Zero}

A more serious issue arises when considering events that have prior
probability equal to 0. In a discrete context this implies that the event in
question is impossible, and so can be ignored, but when using continuous
probability measures this arises quite naturally. It is probably best to think
of a continuous probability as serving as an approximation to something that
is fundamentally discrete as we always measure to finite accuracy. For
example, if each $P_{\theta}$ is absolutely continuous with respect to support
measure $\mu,$ then letting $N_{\delta}(x)$ be a sequence of neighborhoods of
$x$ that converge nicely to $x$ as $\delta\rightarrow0$ (see Rudin (1974)
Chapter 8 and note balls satisfy this) then
\begin{equation}
f_{\theta}(x)=\lim_{\delta\rightarrow0}\frac{P_{\theta}(N_{\delta}(x))}%
{\mu(N_{\delta}(x))} \label{dens1}%
\end{equation}
exists a.e. $\mu$ and serves as a density for $P_{\theta}$ with respect to
$\mu.$ Indeed, when $\mu$ is measuring volume, then (\ref{dens1}) expresses
the fact that $f_{\theta}(x)$ represents the amount of probability per unit
volume at $x.$ In particular, if a version of the density exists which is
continuous and positive at $x,$ then (\ref{dens1}) gives the value of this
density at $x.$ The importance of (\ref{dens1}) is that the value of the
density at $x$ can be thought of as a surrogate for the probability of $x$ via
the approximation $P_{\theta}(N_{\delta}(x))\approx f_{\theta}(x)\mu
(N_{\delta}(x))$ for small $\delta.$ For example, if $\mu$ is volume measure,
then it makes sense to say $x$ is more probable than $y$ when $f_{\theta
}(x)>f_{\theta}(y)$ and the ratio of these quantities reflects the degree.

The natural way to generalize the Bayes factor and relative belief ratio to
the continuous context, where the relevant events have prior content 0, is via
such limits. So let $N_{\delta}(\theta)$ be a sequence of neighborhoods
converging nicely to $\theta$. Then the relative belief ratio at $\theta$ can
be defined as
\[
RB(\theta\,|\,x)=\lim_{\delta\rightarrow0}\frac{\Pi(N_{\delta}(\theta
)\,|\,x)}{\Pi(N_{\delta}(\theta))}%
\]
since the posterior is always absolutely continuous with respect to the prior.
As such, the relative belief ratio is the density of the posterior with
respect to the prior and it exists a.e. $\Pi.$ When the prior has density
$\pi$ and posterior density $\pi(\cdot\,|\,x)$ with respect to support measure
$\upsilon,$ then we obtain
\begin{equation}
RB(\theta\,|\,x)=\frac{\pi(\theta\,|\,x)}{\pi(\theta)}, \label{dens2}%
\end{equation}
the ratio of the densities, whenever $\pi$ is continuous and positive at
$\theta.$ Note too that this definition agrees with the earlier definition
whenever $\Pi(\{\theta\})>0.$ So (\ref{dens2}) is a valid expression of the
principle of evidence in general contexts.

If interest is in the marginal parameter $\psi=\Psi(\theta),$ where
$\Psi:\Theta\rightarrow\Psi,$ with the same notation used for the function and
its range, then the same limiting argument and conditions lead to the relative
belief ratio for $\psi$ given by
\begin{equation}
RB_{\Psi}(\psi\,|\,x)=\frac{\pi_{\Psi}(\psi\,|\,x)}{\pi_{\Psi}(\psi)},
\label{dens3}%
\end{equation}
where $\pi_{\Psi}$ and $\pi_{\Psi}(\cdot\,|\,x)$ are the prior and posterior
densities of $\Psi$ with respect to some common support measure\ $\upsilon
_{\Psi}$ on $\Psi.$ The following result is immediate where $\Pi
(\cdot\,|\,\psi)$ is the conditional prior on $\theta$ given $\Psi
(\theta)=\psi$.\smallskip

\noindent\textbf{Proposition 3. }$RB_{\Psi}(\psi\,|\,x)=E_{\Pi(\cdot
\,|\,\psi)}(RB(\theta\,|\,x)).$\smallskip

\noindent This says that the evidence for $\psi$ is the average of the
evidence for $\theta$ values satisfying $\Psi(\theta)=\psi$ with respect to
the conditional prior. In other words, we get the evidence for $\psi$ by
averaging over the values of the nuisance parameters corresponding to that
value. The Savage-Dickey ratio result gives an alternative expression which is
sometimes useful
\begin{equation}
RB_{\Psi}(\psi\,|\,x)=\frac{m(x\,|\,\psi)}{m(x)}, \label{dens4}%
\end{equation}
where $m(x)=\int_{\Theta}f_{\theta}(x)\,\Pi(d\theta)$ and $m(x\,|\,\psi
)=\int_{\Theta}f_{\theta}(x)\,\Pi(d\theta\,|\,\psi)$ are the prior predictive
and conditional prior predictive densities for $x$ given $\Psi(\theta)=\psi,$
respectively. Note that (\ref{dens4}) is just (\ref{reb1}) but expressed in a
general context.

It is always assumed here that there is enough structure such that $\Pi
(\cdot\,|\,\psi)$ can be defined via a limiting argument so there is no
ambiguity about its definition. This leads to the following general formula
for the density of this probability measure with respect to volume measure on
$\Psi^{-1}\{\psi\},$
\[
\pi(\theta\,|\,\psi)=\frac{\pi(\theta)J_{\Psi}(\theta)}{\pi_{\Psi}(\psi)}%
\]
where $\pi_{\Psi}$ is the prior density of $\Psi,J_{\Psi}(\theta)=(\det
(d\Psi^{t}(\theta)d\Psi(\theta)))^{-1/2}$ and $d\Psi$ is the differential of
$\Psi.$ Note that in the discrete case volume measure is just counting measure
and $J_{\Psi}(\theta)\equiv1.$ The relevant details concerning this
prescription of $\Pi(\cdot\,|\,\psi)$ as a limit can be found in Evans (2015),
Appendix A.

Now consider applying the same approach to obtain a general definition of the
Bayes factor, namely, define
\begin{align*}
BF_{\Psi}(\psi_{0}\,|\,x)  &  =\lim_{\delta\rightarrow0}\left\{  \frac
{\Pi_{\Psi}(N_{\delta}(\psi_{0})\,|\,x)}{\Pi_{\Psi}(N_{\delta}^{c}(\psi
_{0})\,|\,x)}\right\}  \left\{  \frac{\Pi_{\Psi}(N_{\delta}(\psi_{0}))}%
{\Pi_{\Psi}(N_{\delta}^{c}(\psi_{0})}\right\}  ^{-1}\\
&  =RB_{\Psi}(\psi_{0}\,|\,x)\lim_{\delta\rightarrow0}\frac{\Pi_{\Psi
}(N_{\delta}^{c}(\psi_{0}))}{\Pi_{\Psi}(N_{\delta}^{c}(\psi_{0})\,|\,x))}.
\end{align*}
This gives the answer of Proposition 2 when $0<\Pi_{\Psi}(\{\psi_{0}\})<1$ but
when $\Pi_{\Psi}(\{\psi_{0}\})=0$, then the second factor converges to 1. This
establishes the following result.\smallskip

\noindent\textbf{Proposition 4. }When\textbf{ }the prior density $\pi_{\Psi}$
is continuous and positive at $\psi_{0},$ then the limiting Bayes factor and
relative belief ratio agree at $\psi_{0}.$\smallskip

\noindent Therefore, a natural definition of the Bayes factor in the
continuous case is just the relative belief ratio which supports the relative
belief ratio as the more basic measure of evidence generally.

It is not the case, however, that the above limiting argument is the usual
definition of the Bayes factor when $\Pi_{\Psi}(\{\psi_{0}\})=0.$ In this
situation the following approach is usually taken, for example, when assessing
the hypothesis $H_{0}:\Psi(\theta)=\psi_{0}$. For this it is necessary to
specify in addition a probability $p\in(0,1)$ for $H_{0}$ and a prior measure
$\Pi_{\psi_{0}}$\ on $\Psi^{-1}\{\psi_{0}\}.$ Then the overall prior on
$\Theta$ is given by $\Pi_{\psi_{0},p}=p\Pi_{\psi_{0}}+(1-p)\Pi.$ So,
$\Pi_{\psi_{0},p}$ is a mixture of $\Pi_{\psi_{0}}$ with the prior $\Pi,$ and
it assigns the prior mass $p$ to $H_{0}.$ Letting%
\begin{align*}
m_{\psi_{0},p}(x)  &  =\int_{\Theta}f_{\theta}(x)\,\Pi_{p}(d\theta
)=p\int_{\Theta}f_{\theta}(x)\,\Pi_{\psi_{0}}(d\theta)+(1-p)\int_{\Theta
}f_{\theta}(x)\,\Pi(d\theta)\\
&  =pm_{\psi_{0}}(x)+(1-p)m(x),
\end{align*}
the posterior is given by
\[
\Pi_{\psi_{0},p}(A\,|\,x)=\frac{pm_{\psi_{0}}(x)}{m_{\psi_{0},p}(x)}\Pi
_{\psi_{0}}(A\,|\,x)+\frac{(1-p)m(x)}{m_{\psi_{0},p}(x)}\Pi(A\,|\,x).
\]
which implies%
\begin{align*}
BF_{\psi_{0},\Psi}(\psi_{0}\,|\,x)  &  =\frac{\Pi_{p}(\Psi^{-1}\{\psi
_{0}\}\,|\,x)}{\Pi_{p}(\Psi^{-1}\{\psi_{0}\}^{c}\,|\,x)}\frac{1-p}{p}%
=\frac{m_{\psi_{0}}(x)}{m(x)},\\
RB_{\psi_{0},\Psi}(\psi_{0}\,|\,x)  &  =\frac{\Pi_{p}(\Psi^{-1}\{\psi
_{0}\}\,|\,x)}{p}=\frac{m_{\psi_{0}}(x)}{m_{\psi_{0},p}(x)}.
\end{align*}
We have added the subscript $\psi_{0}$ to the notation to denote the fact that
these quantities are now defined using the mixture prior. Notice that $p,$
while necessary for this definition of the Bayes factor, plays no role in its
value so it may seem superfluous to the use of the Bayes factor as a measure
of evidence. Furthermore, we have the following immediate result from the
Savage-Dickey result when $\Pi_{\psi_{0}}=\Pi(\cdot\,|\,\psi_{0})$ since
$m_{\psi_{0}}=m(x\,|\,\psi_{0})$.\smallskip

\noindent\textbf{Proposition 5. }If $\Pi_{\psi_{0}}=\Pi(\cdot\,|\,\psi_{0}),$
the Bayes factor based on the prior $\Pi_{p}$ is the same as the relative
belief ratio based on the prior $\Pi,$ namely, $BF_{\psi_{0},\Psi}(\psi
_{0}\,|\,x)=RB_{\Psi}(\psi_{0}\,|\,x).$\smallskip

\noindent In particular, when $\Psi(\theta)=\theta,$ then Proposition 5 always
holds because there is only one possible prior conditional on $\theta_{0},$
namely, the distribution degenerate at $\theta_{0}.$ Note, however, that when
using the mixture prior the relative belief ratio at $\psi_{0}$ equals
$m_{\psi_{0}}(x)/m_{\psi_{0},p}(x),$ which does depend on $p,$ and so differs
from the Bayes factor.

It is to be noted that, if $\Pi_{\psi_{0}}\neq\Pi(\cdot\,|\,\psi_{0}),$ then
there is a contradiction in the expression of the prior beliefs as $\Pi$
implies the conditional beliefs $\Pi(\cdot\,|\,\psi_{0})$ on $\Psi^{-1}%
\{\psi_{0}\}$ via the limiting argument, while $\Pi_{\psi_{0},p}$ implies the
conditional beliefs $\Pi_{\psi_{0}}.$ In general, if the restriction
$\Pi_{\psi_{0}}=\Pi(\cdot\,|\,\psi_{0})$ is not made, the two measures of
evidence concerning $H_{0}$ can be quite different. For us the avoidance of
this contradiction in prior beliefs is a logical necessity and so Proposition
5 should always apply. This restriction doesn't entirely avoid problems,
however, as even when $\Pi_{\psi_{0}}=$ $\Pi(\cdot\,|\,\psi_{0})$ any
posterior probabilities quoted, say as part of measuring the strength of the
evidence, will be different and this is due to the mixture prior being used
with $p>0.$

\subsection{Inferences}

One of the virtues of the relative belief ratio is that inferences based on it
are now straightforward.\ For estimation there is the relative belief estimate
$\psi(x)=\arg\sup_{\psi}RB_{\Psi}(\psi\,|\,x),$ the value which maximizes the
evidence in favor, and for assessing the hypothesis $H_{0}:\Psi(\theta
)=\psi_{0},$ then the value $RB_{\Psi}(\psi_{0}\,|\,x)$ tells us immediately
whether there is evidence in favor of or against $H_{0}.$ The uncertainties in
these inferences, however, also need to be assessed.\ It is natural to, in
part, assess the uncertainties by measuring the extent to which the inferences
are believed and this is, of necessity, done via quoting posterior probabilities.

For example, the size and posterior content of the \textit{plausible region}
$Pl_{\Psi}(x)=\{\psi:RB_{\Psi}(\psi\,|\,x)>1\},$ the set of values with
evidence in favor, is a natural measure of the accuracy of $\psi(x).$ For the
probability $\Pi(Pl_{\Psi}(x)\,|\,x)$ measures the a posteriori belief that
the true value of $\psi$ is in $Pl_{\Psi}(x).$ For the hypothesis $H_{0}%
:\Psi(\theta)=\psi_{0},$ the posterior probability $\Pi(\{\psi_{0}\}\,|\,x)$
measures how strongly it is believed that $H_{0}$ holds. For example, if
$RB_{\Psi}(\psi_{0}\,|\,x)>1$ and $\Pi(\{\psi_{0}\}\,|\,x)$ is large, then
there is strong evidence in favor of $H_{0},$ etc. One problem with this
measure is that when the prior probability of $\psi_{0}$ is small (e.g., $0$
in the continuous case) we can expect that $\Pi(\{\psi_{0}\}\,|\,x)$ will also
be small (e.g., $0$ in the continuous case). In such contexts it makes sense
to also quote the posterior probability
\begin{equation}
Str_{\Psi,\psi_{0}}(x)=\Pi(RB_{\Psi}(\psi\,|\,x)\leq RB_{\Psi}(\psi
_{0}\,|\,x)\,|\,x) \label{strength}%
\end{equation}
as a measure of the \textit{strength} of the evidence for $H_{0}.$ For when
$RB_{\Psi}(\psi_{0}\,|\,x)>1$ and $Str_{\Psi,\psi_{0}}(x)\approx1,$ then this
says there is strong evidence in favor of $H_{0},$ since there is little
belief that the true value has a larger relative belief ratio while, when
$RB_{\Psi}(\psi_{0}\,|\,x)<1$ and $Str_{\Psi,\psi_{0}}(x)\approx0,$ then this
says there is strong evidence against $H_{0},$ since there is a large belief
that the true value has a larger relative belief ratio. The assessment of the
strength or accuracy of the inferences determined by the evidence can be
carried out in a number of ways as there is no reason to simply quote a single
number. Measuring the strength of the evidence is associated with calibrating
a measure of evidence and this is more extensively discussed in Section 3.2 as
it represents a point of considerable disagreement between what we are
advocating and the common usage of the Bayes factor.

There are a number of good, even optimal properties for relative belief
inferences, see Evans (2015), but one is immediately apparent from
(\ref{dens3}), namely, the inferences are invariant under smooth
reparameterizations since the relevant Jacobians appear in both the prior and
posterior densities and so cancel in the ratio. This invariance is a property
shared with likelihood inferences.

Inferences based on the Bayes factor in the continuous context using the
mixture prior, are focused on the assessment of a single hypothesis
$H_{0}:\Psi(\theta)=\psi_{0}.$\ If we wanted to compare the evidence at
different values $\psi_{1}$ and $\psi_{2},$ then different priors on $\theta$
must be considered. This is not the case for the relative belief ratio and
that is a distinct advantage. Furthermore, in the discrete context the Bayes
factor based on a single prior can be used so this makes the introduction of
the mixture prior, simply to avoid having to deal with a hypothesis that has
prior probability 0, seem very artificial. In any case, as we have argued
here, the use of the mixture prior is not necessary to give a sensible
definition of the Bayes factor.

\subsection{Relationship with Likelihood Theory and Information Theory}

When interest is in the full model parameter $\theta,$ then $RB(\theta
\,|\,x)=f_{\theta}(x)/m(x)$ so it is a member of the equivalence class of
functions $L(\theta\,|\,x)=cf_{\theta}(x),$ for any constant $c>0,$ known as
the likelihood function. Note, however, that the principle of evidence forces
$c=1/m(x)$ so that the relative belief ratio is a valid measure of evidence at
least with cut-off 1. This implies that the relative belief estimate
$\theta(x)$ is the MLE and $Pl(x)$ is a likelihood region. Likelihood theory
doesn't tell us which likelihood region to use, however, while relative belief
does. It is also possible to report a $\gamma$-relative belief credible region
$C_{\gamma}(x)=\{\theta:RB(\theta\,|\,x)\geq c_{\gamma}(x)\}$ where
$c_{\gamma}(x)=\inf_{c}\Pi(\{\theta:RB(\theta\,|\,x)\leq c\}\,|\,x)\geq
1-\gamma,$ provided $\gamma\leq\Pi(Pl(x)\,|\,x).$ The latter constraint
ensures that every value in $C_{\gamma}(x)$ has evidence in its favor as it
doesn't make sense to report values for which we have evidence against as
credible possible choices for the true value of $\theta.$ Note that each
$C_{\gamma}(x)$ is also a likelihood region.

Interestingly all the derivatives with respect to $\theta$ of $\log
RB(\theta\,|\,x)$ agree with the derivatives of $\log L(\theta\,|\,x)$ so the
sampling distributions, given $\theta,$ of these derivatives are the same. For
example, under the usual conditions, it is the case that $E_{\theta}%
(\partial\log RB(\theta\,|\,x)/\partial\theta)=0$ and $E_{\theta}%
(-\partial^{2}\log RB(\theta\,|\,x)/\partial\theta^{2})=I(\theta),$ the Fisher
information. For a marginal parameter $\psi=\Psi(\theta)$ all these results
also hold when the model $\{f_{\theta}:\theta\in\Theta\}$ is replaced by
$\{m(\cdot\,|\,\psi):\psi\in\Psi\}.$ So $RB_{\Psi}(\psi\,|\,x)=m(x\,|\,\psi
)/m(x)$ is in the equivalence class of the integrated likelihood where the
nuisance parameters have been integrated out using $\Pi(\cdot\,|\,\psi),$ see
Berger et al. (1999) for a discussion of integrated likelihood.

The properties for the relative belief ratio as they relate to the likelihood
apply no matter what the prior is. For the Bayes factor for $\theta,$ it is
required that $\Pi(\{\theta\})=\pi(\theta)>0$ and we distinguish two cases. In
the first case suppose that $\pi$ represents a discrete distribution on
$\Theta$ so%
\[
BF(\theta\,|\,x)=\left\{  \frac{f_{\theta}(x)\pi(\theta)}{m(x)-f_{\theta
}(x)\pi(\theta)}\right\}  \left\{  \frac{1-\pi(\theta)}{\pi(\theta)}\right\}
\]
and there is no relationship, as a function of $\theta,$ with the likelihood
alone. A similar comment applies to $BF_{\Psi}(\psi\,|\,x)$ and the integrated
likelihood. For the second case consider the mixture prior $\Pi_{p}%
=p\Pi_{\theta_{0}}+(1-p)\Pi$ where $\Pi_{\theta_{0}}$ is degenerate at
$\theta_{0}$ and $\Pi$ is a continuous prior on $\Theta.$ Then indeed
$BF(\theta_{0}\,|\,x),$ based on $\Pi_{p},$ equals $RB(\theta_{0}\,|\,x)$
based on $\Pi$ but any calibration of this value based on the posterior would
differ because the posteriors differ, see Section 3.2. If we allowed the prior
to vary with $\theta$ so that $BF(\theta\,|\,x)$ is defined for every
$\theta,$ then there is no unique posterior that we can reference for the
calibration. Similar comments apply to the derivatives of the log Bayes factor
and for marginal parameters.

In general, there isn't the same intimate relationship between the Bayes
factor and likelihood as there is between the relative belief ratio and
likelihood. Our view is that such a close relationship is positive for
relative belief as likelihood methodology has proven its value in many
problems. There are contexts where likelihood inference runs into difficulties
but, when these involve repeated sampling properties for the inferences, these
are avoided by relative belief as the inferences and their accuracy
assessments are conditional given the data.

Royall (1997) has developed an evidential approach to inference based on the
likelihood function alone without any involvement of repeated
sampling.\ Certainly, there is a close agreement with relative belief as the
approach is also based on the idea that characterizing statistical evidence is
the key step in developing a\ theory of inference. On the other hand, pure
likelihood suffers from the lack of an unambiguous cut-off for determining
evidence in favor of or against and the dependence on profile likelihoods for
marginal parameters, which in general are not likelihoods in the sense of
being proportional to the probability of something observed. Also, the
complete elimination of repeated sampling from the discussion prevents model
checking and the reliability assessments for the inferences that this
provides. This is not the case for relative belief where such calculations
play a key role, see the discussion of bias in Evans (2015) and Evans and Guo (2021).

Recall that there is nothing to prevent the use of $\log RB_{\Psi}%
(\psi\,|\,x)$ to determine the inferences for $\psi$ as $\log$ is a strictly
increasing function and so the inferences are the same. It is then the case
that
\begin{equation}
E_{\Pi_{\Psi}(\cdot\,|\,x)}(\log RB_{\Psi}(\psi\,|\,x))=KL(\Pi_{\Psi}%
(\cdot\,|\,x),\Pi_{\Psi}), \label{KL}%
\end{equation}
the Kullback-Leibler divergence between the posterior and the prior
distribution of $\psi$. The posterior average value of the log relative belief
ratio is indicating what the average evidence is for values of $\psi,$ at
least if we take $\log RB_{\Psi}(\psi\,|\,x)$ as the evidence measure for
$\psi.$ Since $KL(\Pi_{\Psi}(\cdot\,|\,x),\Pi_{\Psi})\geq0$ this implies that
there is always more evidence in favor of values of $\psi$ than evidence
against unless of course $\Pi_{\Psi}(\cdot\,|\,x)=\Pi_{\Psi}$ and so nothing
has been learned by observing $x$. Also, the mutual information between $\psi$
and $x\,$\ is given by $E(\log RB_{\Psi}(\psi\,|\,x))$ where the expectation
is now with respect to the joint distribution of $(\psi,x).$ It is also
notable that (\ref{KL}) plays a key role in checking for prior-data conflict
as discussed in Nott et al. (2020).

\subsection{Robustness to the Prior}

Now consider the robustness to the prior of the two measures of evidence.
While the results are presented here for the full parameter $\theta,$
corresponding results also hold for the marginal parameter $\psi.$ For the
comparison, the prior is taken to be the same and such that $\Pi
(\{\theta\})=\pi(\theta)>0$ at the point where the sensitivity is being
assessed. Recall that in the continuous case the two measures of evidence
agree when they are defined as limits or when Proposition 5 applies. The
robustness results derived here for the relative belief ratio also apply in
these circumstances. Note that if $T$ is a minimal sufficient statistic for
the model and taking values in $\mathcal{T},$ then $RB(\theta\,|\,x)=RB(\theta
\,|\,T(x))=f_{\theta,T}(T(x))/m_{T}(T(x))$ where $f_{\theta,T}$ is the
sampling density and $m_{T}$ is the prior predictive density of $T.$

\subsubsection{Linear Contamination of the Prior}

Consider linearly contaminated priors
\[
\pi_{\epsilon}(\theta)=\frac{1-\epsilon}{1-\epsilon+\epsilon c_{q}}\pi
(\theta)+\frac{\epsilon c_{q}}{1-\epsilon+\epsilon c_{q}}\frac{q(\theta
)}{c_{q}}%
\]
for some integrable function $q:\Theta\rightarrow\lbrack0,\infty)$ with
normalizing constant $c_{q}$ so $q_{\ast}=q/c_{q}$ is a probability density.
Then $RB(\theta\,|\,T(x))$ is perturbed to $RB_{\epsilon}(\theta
\,|\,T(x))=f_{\theta,T}(T(x))/m_{T,\epsilon}(T(x))$ where
\[
m_{T,\epsilon}(T(x))=\frac{1-\epsilon}{1-\epsilon+\epsilon c_{q}}%
m_{T}(T(x))+\frac{\epsilon c_{q}}{1-\epsilon+\epsilon c_{q}}m_{T}^{(q_{\ast}%
)}(T(x))
\]
and $m_{T}^{(q_{\ast})}$ is the prior predictive density of $T$ obtained using
prior $q_{\ast}.$ The following result follows easily using
\[
\frac{\partial m_{T,\epsilon}(T(x))}{\partial\epsilon}=\frac{c_{q}}{\left(
1-\epsilon+\epsilon c_{q}\right)  ^{2}}\left\{  m_{T}^{(q_{\ast})}%
(T(x))-m_{T}(T(x))\right\}  .
\]

\noindent\textbf{Proposition 6. }Under linear contamination of the prior%
\[
\left.  \frac{\partial\log RB_{\epsilon}(\theta\,|\,T(x))}{\partial\epsilon
}\right\vert _{\epsilon=0}=c_{q}\left(  1-\frac{m_{T}^{(q_{\ast})}%
(T(x))}{m_{T}(T(x))}\right)
\]
which satisfies
\[
E_{m_{T}}\left(  \left.  \frac{\partial\log RB_{\epsilon}(\theta
\,|\,T(x))}{\partial\epsilon}\right\vert _{\epsilon=0}\right)  =0.
\]

\noindent Of some note is that the logarithmic derivative, which measures the
change in a function at a point as a proportion of its value at the point, is
constant in $\theta.\ $This implies that the relative change in the evidence,
as the prior is perturbed, is the same at every value of $\theta.$ Also, no
matter what $q$ is, the a priori expected value of the logarithmic derivative,
at the assumed prior, is 0. Note too that all\ relative belief inferences
remain the same whether we use the relative belief ratio or its logarithm,
e.g., the relative belief estimate and plausible region remain the same. The
ratio $m_{T}^{(q_{\ast})}(T(x))/m_{T}(T(x))$ will be large for some choices of
$q$ when $T(x)$ lies in the tails of $m_{T},$ a situation referred to as
prior-data conflict.\ This emphasizes a point made in Al-Labadi and Evans
(2017), where relative belief inferences are shown to be optimally robust,
among all Bayesian inferences, to the prior on the parameter of interest but
that even these inferences will not be robust if there is prior-data conflict.

It is immediate from Proposition 6 that%
\[
\left.  \frac{\partial RB_{\epsilon}(\theta\,|\,T(x))}{\partial\epsilon
}\right\vert _{\epsilon=0}=c_{q}\left(  1-\frac{m_{T}^{(q_{\ast})}%
(T(x))}{m_{T}(T(x))}\right)  RB(\theta\,|\,T(x)).
\]
Therefore, as $RB(\theta\,|\,T(x))\rightarrow0$, the derivative of the
perturbed relative belief ratio at 0 goes to 0. This indicates that, when
extreme evidence against a value of $\theta$ is obtained, then this finding is
quite robust to the choice of prior. Note that Proposition 14 shows that the
strength of the evidence (\ref{strength}) is bounded above by $RB(\theta
\,|\,T(x))$ and so the strength assessment is also robust to the prior when
there is evidence against, see also Al-Labadi and Evans (2017).

Now consider the Bayes factor at $\theta$ and note that
\[
BF_{\epsilon}(\theta\,|\,x)=\frac{(1-\pi_{\epsilon}(\theta))RB_{\epsilon
}(\theta\,|\,x)}{1-\pi_{\epsilon}(\theta)RB_{\epsilon}(\theta\,|\,x)}.
\]
The Appendix contains the algebra for the following result.\smallskip

\noindent\textbf{Proposition 7. }Under linear contamination of the prior
\begin{align*}
\left.  \frac{\partial\log BF_{\epsilon}(\theta\,|\,T(x))}{\partial\epsilon
}\right\vert _{\epsilon=0}  &  =\left.  \frac{\partial\log RB_{\epsilon
}(\theta\,|\,T(x))}{\partial\epsilon}\right\vert _{\epsilon=0}\left(
1+\frac{\pi(\theta)}{1-\pi(\theta)}BF(\theta\,|\,x)\right)  +\\
&  c_{q}\left(  \frac{\pi(\theta)-q_{\ast}(\theta)}{1-\pi(\theta)}\right)
(1-BF(\theta\,|\,x)).
\end{align*}
So the logarithmic derivative of the Bayes factor at $\epsilon=0$ equals
logarithmic derivative of the relative belief ratio at $\epsilon=0$ plus an
additional term that depends on $\theta.$ The prior expectation of this
quantity is unclear. We do, however, have the following result when interested
in assessing the hypothesis $H_{0}:\theta=\theta_{0}$ where it is common in
that case to take $\pi(\theta_{0})=1/2$ and we suppose that $q_{\ast}%
(\theta_{0})=1/2.$\smallskip

\noindent\textbf{Corollary 1. }When $\pi(\theta_{0})=q_{\ast}(\theta
_{0})=1/2,$ then%
\[
\left.  \frac{\partial\log BF_{\epsilon}(\theta_{0}\,|\,T(x))}{\partial
\epsilon}\right\vert _{\epsilon=0}=\left.  \frac{\partial\log RB_{\epsilon
}(\theta_{0}\,|\,T(x))}{\partial\epsilon}\right\vert _{\epsilon=0}\left(
1+BF(\theta_{0}\,|\,x)\right)  .
\]
so the absolute value of the logarithmic derivative of the Bayes factor at
$\epsilon=0$ is greater than the absolute value of the logarithmic derivative
of the relative belief ratio at $\epsilon=0.$ Furthermore,%
\[
E_{m_{T}}\left(  \left.  \frac{\partial\log BF_{\epsilon}(\theta
\,|\,T(x))}{\partial\epsilon}\right\vert _{\epsilon=0}\right)  =c_{q}\left(
E_{m_{T}^{(q_{\ast})}}\left(  RB(\theta_{0}\,|\,t)\right)  -1\right)  .
\]
Proof: The first statement follows from substitution into the formula and the
fact that $BF(\theta_{0}\,|\,x)>0.$ The second statement follows from
\[
E_{m_{T}}\left(  \left.  \frac{\partial\log BF_{\epsilon}(\theta
\,|\,T(x))}{\partial\epsilon}\right\vert _{\epsilon=0}\right)  =c_{q}\left(
1-\int_{\mathcal{T}}\frac{m_{T}^{(q_{\ast})}(t)}{m_{T}(t)}f_{\theta_{0}%
,T}(t)\,dt\right)
\]
and $RB(\theta_{0}\,|\,t)=f_{\theta_{0},T}(t)/m_{T}(t)$ in this case.
$\blacksquare$\smallskip

\noindent Note that $E_{m_{T}^{(q_{\ast})}}\left(  RB(\theta_{0}%
\,|\,t)\right)  =1$ when $q_{\ast}=\pi$ but otherwise this will generally not
be case so the a priori expected value of the logarithmic derivative of the
Bayes factor at $\epsilon=0$ will not equal 0 in contrast to the situation for
the relative belief ratio.

\subsubsection{Geometric Contamination}

Consider now geometrically contaminated priors
\[
\pi_{\epsilon}(\theta)=\frac{\pi^{1-\epsilon}(\theta)q^{\epsilon}(\theta
)}{\int_{\Theta}\pi^{1-\epsilon}(\theta)q^{\epsilon}(\theta)\,d\theta}%
=\frac{r^{\epsilon}(\theta)\pi(\theta)}{\int_{\Theta}r^{\epsilon}(\theta
)\pi(\theta)\,d\theta}=\frac{r^{\epsilon}(\theta)\pi(\theta)}{E_{\pi
}(r^{\epsilon})}%
\]
for some density function $q:\Theta\rightarrow\lbrack0,\infty)$ such that
$\int_{\Theta}\pi^{1-\epsilon}(\theta)q^{\epsilon}(\theta)\,d\theta<\infty$
and let $r(\theta)=q(\theta)/\pi(\theta).$ For $\theta\in\Theta,$ then
$RB(\theta\,|\,T(x))=f_{\theta,T}(T(x))/m_{T}(T(x))$ is perturbed to
$RB_{\epsilon}(\theta\,|\,x)=f_{\theta,T}(T(x))/m_{\epsilon,T}(T(x))$ where
\[
m_{\epsilon,T}(T(x))=\frac{\int_{\Theta}r^{\epsilon}(\theta)f_{\theta
,T}(T(x))\pi(\theta)\,d\theta}{E_{\pi}(r^{\epsilon})}=m_{T}(T(x))\frac
{E_{\pi(\cdot\,|\,x)}(r^{\epsilon})}{E_{\pi}(r^{\epsilon})}.
\]
Therefore, assuming the conditions of the dominated derivative theorem,%
\begin{align*}
\frac{1}{m_{\epsilon,T}(T(x))}\frac{\partial m_{\epsilon,T}(T(x))}%
{\partial\epsilon}  &  =\frac{E_{\pi(\cdot\,|\,x)}(r^{\epsilon}\log r)}%
{E_{\pi}(r^{\epsilon})}-\frac{E_{\pi(\cdot\,|\,x)}(r^{\epsilon})}{E_{\pi}%
^{2}(r^{\epsilon})}E_{\pi}(r^{\epsilon}\log r)\\
&  \rightarrow E_{\pi(\cdot\,|\,x)}(\log r)-E_{\pi}(\log r)\text{ as }%
\epsilon\rightarrow0.
\end{align*}
This implies the following result.\smallskip

\noindent\textbf{Proposition 8. }Under geometric contamination of the prior
\[
\left.  \frac{\partial\log RB_{\epsilon}(\theta\,|\,x)}{\partial\epsilon
}\right\vert _{\epsilon=0}=E_{\pi}(\log r)-E_{\pi(\cdot\,|\,x)}(\log r)
\]
which satisfies
\[
E_{m_{T}}\left(  \left.  \frac{\partial\log RB_{\epsilon}(\theta
\,|\,x)}{\partial\epsilon}\right\vert _{\epsilon=0}\right)  =0.
\]

\noindent Again, this derivative does not depend on $\theta$ and has prior
expectation equal to 0.

The following gives the result for the Bayes factor with proof in the
Appendix.\smallskip

\noindent\textbf{Proposition 9. }Under geometric contamination of the prior
\begin{align*}
\left.  \frac{\partial\log BF_{\epsilon}(\theta\,|\,x)}{\partial\epsilon
}\right\vert _{\epsilon=0}  &  =-\left.  \frac{\partial\log RB_{\epsilon
}(\theta\,|\,x)}{\partial\epsilon}\right\vert _{\epsilon=0}\left(  1+\frac
{\pi(\theta)}{1-\pi(\theta)}BF(\theta\,|\,x)\right)  +\\
&  (E_{\pi}(\log r)-\log r(\theta))\left(  \frac{\pi(\theta)}{1-\pi(\theta
)}\right)  \left(  1-BF(\theta\,|\,x)\right)  .
\end{align*}

\noindent This clearly depends on $\theta.$ In the most common usage of the
Bayes factor the prior odds in favor of $\theta_{0}$ are taken to be 1 and the
following obtains.\smallskip

\noindent\textbf{Corollary 2. }If $\pi(\theta_{0})=q(\theta_{0})=1/2,$ then%
\begin{align*}
&  \left.  \frac{\partial\log BF_{\epsilon}(\theta_{0}\,|\,x)}{\partial
\epsilon}\right\vert _{\epsilon=0}\\
&  =-\left.  \frac{\partial\log RB_{\epsilon}(\theta_{0}\,|\,x)}%
{\partial\epsilon}\right\vert _{\epsilon=0}\left(  1+BF(\theta_{0}%
\,|\,x)\right)  +E_{\pi}(\log r)\left(  1-BF(\theta_{0}\,|\,x)\right)  .
\end{align*}

\noindent There doesn't appear to be any further simplification here although
the first term has absolute value bigger than the corresponding term for the
log relative belief ratio.

\subsection{Information Inconsistency}

The following example is also discussed in Evans (2015) and illustrates this
issue. Suppose $x=(x_{1},\ldots,x_{n})$ is a sample from a $N(\mu,\sigma^{2})$
distribution, where $\mu\in R^{1},\sigma^{2}>0$ are unknown, so $T(x)=(\bar
{x},s^{2}),$ where $\sqrt{n}\bar{x}\sim N(\mu,1)\ $independently of
$s^{2}=(x-\bar{x}1)^{\prime}(x-\bar{x}1)\sim\sigma^{2}$chi-squared$(n-1),$ is
a minimal sufficient statistic.\ Consider the prior $\Pi_{\mu_{0},p}$ where,
with $\delta_{\mu_{0}}$ the probability measure on $%
\mathbb{R}
^{1}$ degenerate at $\mu_{0},$%
\[
\mu\,|\,\sigma^{2}\sim p\delta_{\mu_{0}}+(1-p)N(\mu_{0},\tau_{0}^{2}\sigma
^{2}),\quad1/\sigma^{2}\sim\text{ gamma}_{\text{rate}}(\alpha_{0},\beta_{0}),
\]
so $p,\tau_{0}^{2},\alpha_{0},\beta_{0}$ are hyperparameters. Suppose
$\Psi(\mu)=\mu,$ and we want to assess the hypothesis $H_{0}:\mu=\mu_{0}.$
Some straight-forward computations give%
\[
BF_{\mu_{0},\Psi}(\mu_{0}\,|\,x)=(1+n\tau_{0}^{2})^{1/2}\left[  \frac
{(1+n\tau_{0}^{2})^{-1}n(\bar{x}-\mu_{0})^{2}+s^{2}+2\beta_{0}}{n(\bar{x}%
-\mu_{0})^{2}+s^{2}+2\beta_{0}}\right]  ^{n/2+\alpha_{0}}.
\]
Notice that as $\sqrt{n}\left\vert \bar{x}-\mu_{0}\right\vert \rightarrow
\infty$ with $s^{2}$ fixed, then $BF_{\mu_{0},\Psi}(\mu_{0}\,|\,x)\rightarrow
(1+n\tau_{0}^{2})^{-n/2-\alpha_{0}+1/2}\neq0$ which is counter-intuitive as
one might expect that, as the data become more and more divergent from the
hypothesized value, the Bayes factor should converge to 0 and give categorical
evidence against $H_{0}$. This phenomenon is not uncommon and it is referred
to as \textit{information inconsistency}.

Now examine the relative belief ratio for this problem based on the continuous
prior $\Pi$ where $\mu\,|\,\sigma^{2}\sim N(\mu_{0},\tau_{0}^{2}\sigma
^{2}),1/\sigma^{2}\sim$ gamma$_{\text{rate}}(\alpha_{0},\beta_{0}).$ A
straightforward calculation gives%
\begin{align*}
&  RB_{\Psi}(\mu_{0}\,|\,x)\\
&  =\left\{  \frac{\Gamma\left(  \frac{n+2\alpha_{0}}{2}\right)  }%
{\Gamma\left(  \frac{1}{2}\right)  \Gamma\left(  \frac{n-1+2\alpha_{0}}%
{2}\right)  }/\frac{\Gamma\left(  \frac{1+2\alpha_{0}}{2}\right)  }%
{\Gamma\left(  \frac{1}{2}\right)  \Gamma\left(  \frac{2\alpha_{0}}{2}\right)
}\right\}  \left\{  \frac{\sqrt{\tau_{0}^{2}\beta_{0}/\alpha_{0}}}%
{\sqrt{(n+1/\tau_{0}^{2})^{-1}\beta(x)/(n+2\alpha_{0})}}\right\}  \times\\
&  \left\{  (1+(n+2\alpha_{0})(\mu_{0}-\mu(x))^{2}/(n+1/\tau_{0}^{2}%
)^{-1}\beta(x))^{-\frac{n+1+2\alpha_{0}}{2}}\right\}
\end{align*}
where%
\begin{align*}
\mu(x)  &  =(n+1/\tau_{0}^{2})^{-1}(n\bar{x}+\mu_{0}/\tau_{0}^{2}%
)=(n+1/\tau_{0}^{2})^{-1}n(\bar{x}-\mu_{0})+\mu_{0},\\
\beta(x)  &  =(n+1/\tau_{0}^{2})^{-1}[n(\bar{x}-\mu_{0})^{2}/\tau_{0}%
^{2}]+s^{2}+2\beta_{0}.
\end{align*}
It is then seen that $RB_{\Psi}(\mu_{0}\,|\,x)\rightarrow0$ as $\sqrt{n}%
|\bar{x}-\mu_{0}|\rightarrow\infty$ with $s^{2}$ constant. So, the relative
belief ratio avoids information inconsistency here.

This particular example generalizes to the regression context so this is not
an isolated example. Mulder et al. (2017) considers this phenomenon and notes
that information inconsistency can be avoided with certain priors. They do
not, however, provide a reason for the inconsistency. So it would seem that we
are not allowed to choose just any prior if it is desirable to avoid this
phenomenon. Since the prior $\Pi$ is conjugate this seems anomalous but what
has been shown here is that the problem is avoided if we simply use the
continuous prior together with the relative belief ratio as opposed to the
mixture prior and the Bayes factor.

\section{Criticisms of these Measures of Evidence}

Criticisms have been made in literature concerning Bayes factors and typically
the same comments can be applied to the relative belief ratio as well. In our
view, these criticisms are based on an incomplete view of what it means to
measure statistical evidence as is now discussed.

\subsection{Incoherence}

This can be demonstrated using a probability model $(\Omega,\mathcal{A},P)$
and events $A,B,C$ all having positive probability with $A\subset
B,P(B\backslash A)>0$ and suppose $C$ is observed to be true. The incoherence
supposedly arises because it can happen that $RB(A\,|\,C)>RB(B\,|\,C),$ and
similarly for the Bayes factor, so these measures are not monotonic and some
argue that they should be. Certainly, probability is monotonic, which is
appropriate as it is measuring degree of belief. Measures of evidence,
however, are based on change in belief because, whether or not the occurrence
of the event $C$ is evidence in favor of or against an event being true,
depends on whether the conditional belief given $C$ is higher or lower than
the initial belief. There is the following result which helps to understand
this phenomenon.\smallskip

\noindent\textbf{Proposition 10. }When $A\subset B$ and all the events
$A,B\backslash A$ and $C$ have positive probability, then
$RB(A\,|\,C)>RB(B\,|\,C)$ iff $RB(A\,|\,C)>RB(B\backslash A\,|\,C).$

\noindent Proof: Now
\begin{align*}
RB(A\,|\,C)  &  =\frac{P(A\,|\,C)}{P(A)}>RB(B\,|\,C)=\frac{P(B\,|\,C)}%
{P(B)}\text{ iff }\\
\frac{P(B)}{P(A)}  &  =\frac{P(B\backslash A)}{P(A)}+1>\frac{P(B\backslash
A\,|\,C)}{P(A\,|\,C)}+1
\end{align*}
which gives the result. $\blacksquare$\smallskip

\noindent A simple example makes it clear why this result makes
sense.\smallskip

\noindent\textbf{Example 3.1.1 }

Suppose $\Omega$ is the population of a country where $B$ is the set of adults
over 20 years of age and $A\subset B$ is the set of university graduates over
20. An individual $\omega$ is randomly selected and it is determined that
$\omega\in C,$ the set of those in favor of vaccination against COVID-19.
Certainly we might expect there to be evidence in favor of $\omega\in A$ and
of $\omega\in B$ as it is certainly reasonable that the proportion of
university graduates over 20 in favor of vaccination, within the subpopulation
of those who believe in vaccination, is larger than proportion of university
graduates within the population of those over 20, namely, $P(A\,|\,C)>P(A),$
and this implies that observing $\omega\in C$ is evidence in favor of
$\omega\in A.$ A similar comment might also apply to $B.$ It would not be at
all surprising, however, that $RB(A\,|\,C)>RB(B\backslash A\,|\,C)$ since
$B\backslash A$ are those over 20 with less education and it might even be the
case that $P(B\backslash A\,|\,C)<P(B\backslash A)$ so there is evidence
against $B.$ $\blacksquare$\smallskip

Perhaps it is even more surprising that the occurrence of $C$ may lead to
evidence in favor of $A\subset B$ but evidence against $B.$ The following
result, established in Evans (2015), indicates when this will occur.\smallskip

\noindent\textbf{Proposition 11. }Under the\textbf{ }conditions of Proposition
12\textbf{, }if\textbf{ }$RB(A\,|\,C)>1,$ then $RB(B\,|\,C)<1$ iff
$RB(B\backslash A\,|\,C)<1$ and
\begin{equation}
P(A\,|\,B)<\frac{1-RB(B\backslash A\,|\,C)}{RB(A\,|\,C)-RB(B\backslash
A\,|\,C)}. \label{lem121}%
\end{equation}

\noindent Proof: Under the conditions of Proposition 12%
\[
RB(B\,|\,C)=RB(A\,|\,C)P(A\,|\,B)+RB(B\backslash A\,|\,C)P(B\backslash
A\,|\,B)
\]
and then $RB(B\,|\,C)<1$ forces $RB(B\backslash A\,|\,C)<1$ since
$RB(A\,|\,C)>1$ and $P(A\,|\,B)>0.$ Then $RB(A\,|\,C)P(A\,|\,B)+RB(B\backslash
A\,|\,C)(1-P(A\,|\,B))<1$ iff
\[
(RB(A\,|\,C)-RB(B\backslash A\,|\,C))P(A\,|\,B)<1-RB(B\backslash A\,|\,C)
\]
which gives (\ref{lem121}). Now suppose $RB(B\backslash A\,|\,C)<1$ and
(\ref{lem121}) holds. Then from (\ref{lem121})%

\[
RB(B\,|\,C)=(RB(A\,|\,C)-RB(B\backslash A\,|\,C))P(A\,|\,B)+RB(B\backslash
A\,|\,C)<1.\text{ }\blacksquare
\]

\noindent So this apparent anomaly occurs whenever $P(A\,|\,B)$ is very small
or $RB(A\,|\,C)\approx1.$ In terms of Example 3.1, this will arise in a
country where there is evidence against $\omega\in B\backslash A$ and either
the set of university graduates over 20 is a very small proportion of those
over 20 or those over 20 who are university educated virtually all believe in vaccination.

For us there is no incoherency in the relative belief ratio as a measure of
evidence.\ The discussion here shows that thinking about evidence is more
subtle than perhaps first expected.

\subsection{Measuring the Strength of Evidence: Calibration}

Some criticisms of Bayes factors have been based on comparing Bayes factors
from quite different contexts. For example, Sarafoglou et al. (2022) discuss a
situation where two Bayes factors are compared for essentially the same
situation but with different sample sizes and respond to a criticism that is
similar to the monotonicity issue discussed in Section 3.1. Also, consider the
example of Section 2.6 and notice that as $\tau_{0}^{2}\rightarrow\infty,$
then $BF_{\mu_{0},\Psi}(\mu_{0}\,|\,x)\rightarrow\infty$ and $RB_{\Psi}%
(\mu_{0}\,|\,x)\rightarrow\infty.$ If we take the value of the measure of
evidence as also measuring the strength of the evidence, then this would seem
to indicate that it is possible to obtain overwhelming evidence in favor of
$H_{0}$ by taking the prior to be sufficiently diffuse and this is true no
matter how large $t=\sqrt{n}|\bar{x}-\mu_{0}|/(s/\sqrt{n-1})$ is. This is
generally referred to as the Jeffreys-Lindley paradox because a frequentist
would find evidence against $H_{0}$ when $t$ is large but a sufficiently
diffuse Bayesian would not.

As with estimation problems, where the accuracy of an estimate also has to be
assessed via some additional measure, it is our view that in hypothesis
assessment problems the strength of evidence, whether in favor of or against,
must also be assessed and assessed separately from the value of the measure of
evidence. We refer to this as \textit{calibrating} the measure of evidence. As
discussed in Kass and Raftery (1995) there is scale, due to Jeffreys, see
Table \ref{tab1}, that is to be used to assess strength of evidence given by
the Bayes factor but there doesn't appear to be a good supporting argument for
this scale.
\begin{table}[tbp] \centering
\begin{tabular}
[c]{|l|l|}\hline
$BF$ & Strength\\\hline
$1$ to $10^{1/2}$ & Barely worth mentioning\\\hline
$10^{1/2}$ to $10$ & Substantial\\\hline
$10$ to $10^{3/2}$ & Strong\\\hline
$10^{3/2}$ to $10^{2}$ & Very Strong\\\hline%
$>$
$10^{2}$ & Decisive\\\hline
\end{tabular}
\caption{Jeffreys' Bayes factor scale for measuring strength of the evidence in favor (the strength of evidence against is measured by the reciprocals),}\label{tab1}%
\end{table}%

Various authors have suggested modifications of Jeffreys' scale but none of
these deal with the essential question: is there a universal scale on which
statistical evidence can be measured such that, for two very different
contexts $BF_{1}=BF_{2}$ means, not only that both contexts have found
evidence against or evidence in favor, but also that the strength of this
evidence is the same? Certainly the phenomenon of the Jeffreys-Lindley paradox
suggests otherwise but the following very simple example makes the answer
quite clear that the answer is definitively no.\smallskip

\noindent\textbf{Example 3.2.1} \textit{Prosecutor's Fallacy}

Suppose that a uniform probability distribution is placed on a population of
size $N$ and that some trait left at a crime scene is shared by $m\ll N$
members of the population. An individual is randomly selected from the
population and it is found that they possess the trait. A prosecutor cites the
rarity of the trait as strong evidence of guilt of this individual. It is to
be noted, however, that the posterior probability of guilt is $P($%
\textquotedblleft guilty\textquotedblright$\,|\,$\textquotedblleft has
trait\textquotedblright$)=1/m$\ and this could be very small, e.g., when
$m=10^{3}.$ In essence it seems that the prosecutor has flipped the antecedent
and the consequent in drawing their conclusion as indeed $P($\textquotedblleft
has trait\textquotedblright$\,|\,$\textquotedblleft
guilty$\text{\textquotedblright})=1.$

It has to be acknowledged, however, that the possession of the trait is
evidence of guilt so the prosecutor isn't completely wrong. Indeed,%

\begin{align*}
RB(\text{\textquotedblleft guilty\textquotedblright}%
\,|\,\text{\textquotedblleft has trait\textquotedblright})  &  =\frac
{P(\text{\textquotedblleft guilty\textquotedblright}%
\,|\,\text{\textquotedblleft has trait\textquotedblright})}%
{P(\text{\textquotedblleft guilty\textquotedblright)}}=\frac{N}{m}\gg1\\
BF(\text{\textquotedblleft guilty\textquotedblright}%
\,|\,\text{\textquotedblleft has trait\textquotedblright})  &  =\frac
{RB(\text{\textquotedblleft guilty\textquotedblright}%
\,|\,\text{\textquotedblleft has trait\textquotedblright})}%
{RB(\text{\textquotedblleft innocent\textquotedblright}%
\,|\,\text{\textquotedblleft has trait\textquotedblright})}=\frac{N-1}{m-1}%
\gg1.
\end{align*}
Both measures indicate evidence in favor of guilt as they should and the Bayes
factor is bigger than the relative belief ratio as it always is when there is
evidence in favor. But we now ask how strongly do we believe what the evidence
says and this is measured by $P($\textquotedblleft guilty\textquotedblright%
$\,|\,$\textquotedblleft has trait\textquotedblright$)=1/m.$ If $m=1,$ the
evidence is categorical but if $m=10^{3}$ the evidence can only be considered
as being very weak. Certainly, most would not want to convict someone based on
such a weak belief in what the evidence indicates.

Now consider various values of $N$ and suppose $m=10^{3}.$ If $N=10^{5},$ then
$BF=100.01$ which according to Jeffreys scale is definitive evidence in favor
of guilt. If $N=10^{6},$ then $BF=1001.00$ which is even more definitive
according to this scale. Clearly, one can construct examples where the Bayes
factor can be as large as one chooses but this does not make the evidence any
more definitive.

Of course, the same criticism could be leveled at the relative belief ratio
\textit{if} we were advocating using its size to measure its strength, but
that is not the case. It is most natural to measure the strength of the
evidence by a posterior probability which is measuring how strongly we believe
what the evidence says and in this simple binary situation of guilt versus
innocence, the posterior probability $P($\textquotedblleft
guilty\textquotedblright$\,|\,$\textquotedblleft has trait\textquotedblright%
$)$ accomplishes this.

One general criticism that can be raised concerning how we propose to measure
the strength of evidence is as follows: what values of the posterior
probability constitute strong and weak evidence? The answer to this lies in
the application and what the ultimate outcome of the inference will be. For
example, in the context of this example a posterior probability of $1/10^{3}$
seems far too small to conclude guilt and convict even though there is
evidence of guilt. On the other hand, suppose the same numbers apply in a
situation where individuals were being tested to see if they were carriers of
a highly infectious disease and a positive test was very weak evidence in
favor of that being the case. If the outcome was to quarantine an individual
for a few days, then even very weak evidence could justify this. There are
indeed situations where costs come into play in determining actions/decisions
but this has nothing to do with what the evidence says and how strong it is.
$\blacksquare$\smallskip

In general, there does not appear to be a universal scale on which evidence
and its strength are to be measured and so we advocate that such an assessment
be made based on the particular context. This means that, for example, a
relative belief ratio of 2 in different contexts, certainly means that
evidence in favor has been found, but the strength of this evidence may be
quite different.

Exactly how the strength of the evidence is to be measured is not entirely
clear and there is no reason to suppose that only one measure be used.\ But
consider a typical problem where the hypothesis $H_{0}:\Psi(\theta)=\psi_{0}$
is to be assessed. Then, presuming the values $\psi\in\Psi$ all correspond to
some real world object with a common interpretation, it makes sense to compare
$RB_{\Psi}(\psi_{0}\,|\,x)$ to the other values $RB_{\Psi}(\psi\,|\,x)$ for
$\psi\in\Psi$ and that is what (\ref{strength}) does with different
interpretations depending on whether $RB_{\Psi}(\psi_{0}\,|\,x)$ indicates
evidence in favor of or against $H_{0}.$ Of some importance here is that, even
within a given context, the comparison be made among like objects so we aren't
comparing apples with oranges. For example, if $\Psi=I_{A},$ the indicator of
the event $A\subset\Theta$ with positive prior content, and with $\psi_{0}=1,$
then quoting the posterior probability $\Pi(A\,|\,x)$ is often a more
appropriate measure of the strength of the evidence than (\ref{strength})
unless there is a partition of $\Theta$ into similar subsets with $A$ a
member. When such a partition doesn't exist, if $RB_{\Psi}(A\,|\,x)>1$ and
$\Pi(A\,|\,x)$ is large, then there is strong evidence in favor of $A$ being
true, etc. A general alternative for measuring the strength of evidence is to
quote $Pl_{\Psi}(x)$ and its posterior content. For if $RB_{\Psi}(\psi
_{0}\,|\,x)>1,$ then $\psi_{0}\in Pl_{\Psi}(x)$ and if $Pl_{\Psi}(x)$ has a
small size and $\Pi(Pl_{\Psi}(x)\,|\,x)$ is large, this can be interpreted as
strong evidence in favor of $H_{0},\ $etc. It is worth noting too that it also
possible to calibrate a Bayes factor in a similar way in the purely discrete
context where a prior does not require modification. The comparative
simplicity of the relative belief ratio to measure evidence and
(\ref{strength}) to measure the strength of evidence makes this approach preferable.

Several inequalities can provide further insight into the size issue when
considering the relative belief ratio as in the following result.\smallskip

\noindent\textbf{Proposition 12. }The following posterior inequalities hold
\begin{align}
&  1-\min\left\{  \frac{E_{\Pi_{\Psi}(\cdot\,|\,x)}(RB_{\Psi}(\psi
\,|\,x))}{RB_{\Psi}(\psi_{0}\,|\,x)},\frac{KL(\Pi_{\Psi}(\cdot\,|\,x),\Pi
_{\Psi})}{\log RB_{\Psi}(\psi_{0}\,|\,x)}\right\} \nonumber\\
&  \leq Str_{\Psi,\psi_{0}}(x)\leq RB_{\Psi}(\psi_{0}\,|\,x) \label{evineq}%
\end{align}
Proof: Using Markov's inequality we have
\begin{align*}
&  \Pi_{\Psi}(RB_{\Psi}(\psi\,|\,x)\leq RB_{\Psi}(\psi_{0}\,|\,x)\,|\,x)=\Pi
_{\Psi}(RB_{\Psi}^{-1}(\psi\,|\,x)\geq RB_{\Psi}^{-1}(\psi_{0}\,|\,x)\,|\,x)\\
&  \leq\frac{E_{\Pi_{\Psi}(\cdot\,|\,x)}(RB_{\Psi}^{-1}(\psi\,|\,x))}%
{RB_{\Psi}^{-1}(\psi_{0}\,|\,x)}=RB_{\Psi}(\psi_{0}\,|\,x),
\end{align*}
since $E_{\Pi(\cdot\,|\,x)}(RB_{\Psi}^{-1}(\psi\,|\,x))=1.$ For the left-hand
side we have%
\begin{align*}
&  \Pi_{\Psi}(RB_{\Psi}(\psi\,|\,x)\leq RB_{\Psi}(\psi_{0}\,|\,x)\,|\,x)=1-\Pi
_{\Psi}(RB_{\Psi}(\psi\,|\,x)>RB_{\Psi}(\psi_{0}\,|\,x)\,|\,x)\\
&  \geq1-\frac{E_{\Pi_{\Psi}(\cdot\,|\,x)}(RB_{\Psi}(\psi\,|\,x))}{RB_{\Psi
}(\psi_{0}\,|\,x)}%
\end{align*}
and
\begin{align*}
&  \Pi_{\Psi}(RB_{\Psi}(\psi\,|\,x)\leq RB_{\Psi}(\psi_{0}\,|\,x)\,|\,x)=1-\Pi
_{\Psi}(\log RB_{\Psi}(\psi\,|\,x)>\log RB_{\Psi}(\psi_{0}\,|\,x)\,|\,x)\\
&  \geq1-\frac{KL(\Pi_{\Psi}(\cdot\,|\,x),\Pi_{\Psi})}{\log RB_{\Psi}(\psi
_{0}\,|\,x)\,|\,x)}.\text{ }\blacksquare
\end{align*}
The right hand side of (\ref{evineq}) indicates that, if there is evidence
against $\psi_{0}$ and $RB_{\Psi}(\psi_{0}\,|\,x)<<1,$ then this corresponds
to strong evidence against $\psi_{0}$ being the true value and so there is no
need to compute the strength. The left-hand side indicates that, if there is
evidence in favor of $\psi_{0}$ and $RB_{\Psi}(\psi_{0}\,|\,x)$ is big
relative to the size of $E_{\Pi_{\Psi}(\cdot\,|\,x)}(RB_{\Psi}(\psi\,|\,x)),$
or $\log RB_{\Psi}(\psi_{0}\,|\,x)\,|\,x)$ is big relative to the size of
$KL(\Pi_{\Psi}(\cdot\,|\,x),\Pi_{\Psi})$, then there is strong evidence in
favor of $\psi_{0}.$

Interestingly, there is a clear asymmetry between the evidence against and the
evidence in favor situations. For a very small value of $RB_{\Psi}(\psi
_{0}\,|\,x)$ is always strong evidence against $\psi_{0}$ no matter what the
other values of $RB_{\Psi}(\cdot\,|\,x)$ may be. For evidence in favor,
however, the value of $RB_{\Psi}(\psi_{0}\,|\,x),$ no matter how big or small,
always has to be calibrated against the other values of $RB_{\Psi}%
(\cdot\,|\,x)$. For example, a very large value of $RB_{\Psi}(\psi_{0}\,|\,x)$
does not mean strong evidence in favor, as exemplified by Example 3.2.1, while
even a value of $RB_{\Psi}(\psi_{0}\,|\,x)$ just bigger than 1 can be strong
evidence in favor if $KL(\Pi_{\Psi}(\cdot\,|\,x),\Pi_{\Psi})$ is small,
namely, if the data have not changed beliefs by much. This latter situation
raises a general issue concerning the reliability of Bayesian inferences which
is a topic addressed in Evans (2015) and Evans and Guo (2021) and is
essentially a design issue ensuring that a sufficient amount of data has been
obtained to trust the inferences.

\section{Conclusions}

This paper has been concerned with comparing the Bayes factor and the relative
belief ratio as measures of evidence. The overall conclusion is that the
relative belief ratio is much more suited to this task. There does not appear
to be any reason why the Bayes factor should be preferred. The relative belief
ratio is simpler and has nicer properties. The fact that current usage of the
Bayes factor for hypothesis assessment requires the modification of a
perfectly good prior in the continuous case to a mixture prior with a discrete
component, is in itself an indication that there is something amiss. This does
not mean that such a mixture prior is inappropriate in certain contexts, it is
the \textit{necessity} to make this modification so that the Bayes factor can
be used, that is objectionable. Even in contexts where such a mixture prior is
used, the relative belief ratio has better robustness to the prior properties
as demonstrated here. Another nice feature of the relative belief ratio is
that it can be used for hypothesis assessment, estimation and prediction while
this is currently not the case for the Bayes factor. This makes sense as any
proper characterization of statistical evidence must be applicable to all
statistical inference problems. Also, information inconsistency associated
with the Bayes factor based on the mixture prior can be avoided by using the
relative belief ratio with the associated continuous prior. In effect, as the
Bayesian approach to inference requires, with relative belief we are free to
use any prior and are not constrained by any mathematical phenomenon such as
the need to have a discrete mass at a point or the need to avoid information inconsistency.

No discussion of the measurement of statistical evidence in the frequentist
context has been provided here. This is because there is nothing as simple as
the principle of evidence that can be used to characterize statistical
evidence without a prior. So this remains an open problem in that context.

\section*{Appendix}

\noindent\textbf{Proof of Proposition 7} Using $\partial\pi_{\epsilon}%
(\theta)/\partial\epsilon=c_{q}(q_{\ast}(\theta)-\pi(\theta))/(1-\epsilon
+\epsilon c_{q})^{2},$%
\begin{align}
\frac{\partial\log BF_{\epsilon}(\theta\,|\,x)}{\partial\epsilon}  &
=\frac{-\frac{\partial\pi_{\epsilon}(\theta)}{\partial\epsilon}}%
{1-\pi_{\epsilon}(\theta_{0})}+\frac{\partial\log RB_{\epsilon}(\theta
\,|\,x)}{\partial\epsilon}+\nonumber\\
&  \frac{RB_{\epsilon}(\theta\,|\,x)\frac{\partial\pi_{\epsilon}(\theta
)}{\partial\epsilon}+\pi_{\epsilon}(\theta)\frac{\partial RB_{\epsilon}%
(\theta\,|\,x)}{\partial\epsilon}}{1-\pi_{\epsilon}(\theta)RB_{\epsilon
}(\theta\,|\,x)}, \label{logderiv}%
\end{align}
which implies%
\begin{align*}
&  \left.  \frac{\partial\log BF_{\epsilon}(\theta_{0}\,|\,T(x))}%
{\partial\epsilon}\right\vert _{\epsilon=0}=c_{q}\left(  \frac{\pi
(\theta)-q_{\ast}(\theta)}{1-\pi(\theta)}\right)  +c_{q}\left(  1-\frac
{m^{(q_{\ast})}(x)}{m(x)}\right)  +\\
&  c_{q}\frac{(1-\pi(\theta))RB(\theta\,|\,x)}{1-\pi(\theta)RB(\theta
\,|\,x)}\left(  \frac{q_{\ast}(\theta)-\pi(\theta)}{1-\pi(\theta)}+\frac
{\pi(\theta)}{1-\pi(\theta)}\left(  1-\frac{m^{(q_{\ast})}(x)}{m(x)}\right)
\right) \\
&  =\left.  \frac{\partial\log RB_{\epsilon}(\theta_{0}\,|\,T(x))}%
{\partial\epsilon}\right\vert _{\epsilon=0}\left(  1+\frac{\pi(\theta)}%
{1-\pi(\theta)}BF(\theta\,|\,x)\right)  +\\
&  c_{q}\left(  \frac{\pi(\theta)-q_{\ast}(\theta)}{1-\pi(\theta)}\right)
(1-BF(\theta\,|\,x)).
\end{align*}

\noindent\textbf{Proof of Proposition 9 }Using $\frac{\partial\pi_{\epsilon
}(\theta)}{\partial\epsilon}=\frac{(\log r(\theta))r^{\epsilon}(\theta
)\pi(\theta)}{E_{\pi}(r^{\epsilon})}-\frac{r^{\epsilon}(\theta)\pi(\theta
)}{E_{\pi}^{2}(r^{\epsilon})}E_{\pi}(r^{\epsilon}\log r)$\newline%
$\rightarrow(\log r(\theta)-E_{\pi}(\log r))\pi(\theta)$ as $\epsilon
\rightarrow0$ and substituting the relevant expressions into (\ref{logderiv})
gives%
\begin{align*}
&  \left.  \frac{\partial\log BF_{\epsilon}(\theta\,|\,x)}{\partial\epsilon
}\right\vert _{\epsilon=0}\\
&  =-(\log r(\theta)-E_{\pi}(\log r))\frac{\pi(\theta)}{1-\pi(\theta)}%
+(E_{\pi(\cdot\,|\,x)}(\log r)-E_{\pi}(\log r))+\\
&  \frac{RB(\theta\,|\,x)(\log r(\theta)-E_{\pi}(\log r))\pi(\theta
)+\pi(\theta)(E_{\pi(\cdot\,|\,x)}(\log r)-E_{\pi}(\log r))RB(\theta
\,|\,x)}{1-\pi(\theta)RB(\theta\,|\,x)}\\
&  =-(\log r(\theta)-E_{\pi}(\log r))\left(  \frac{\pi(\theta)}{1-\pi(\theta
)}-\frac{RB(\theta\,|\,x)\pi(\theta)}{1-\pi(\theta)RB(\theta\,|\,x)}\right)
+\\
&  (E_{\pi(\cdot\,|\,x)}(\log r)-E_{\pi}(\log r))\left(  1+\frac{\pi
(\theta)RB(\theta\,|\,x)}{1-\pi(\theta)RB(\theta\,|\,x)}\right) \\
&  =(E_{\pi}(\log r)-\log r(\theta))\left(  \frac{\pi(\theta)}{1-\pi(\theta
)}\right)  \left(  1-BF(\theta\,|\,x)\right)  -\\
&  (E_{\pi}(\log r)-E_{\pi(\cdot\,|\,x)}(\log r))\left(  1+\frac{\pi(\theta
)}{1-\pi(\theta)}BF(\theta\,|\,x)\right)  .
\end{align*}

\section{References}

\noindent Achinstein, P. (2001) The Book of Evidence. Oxford University
Press.\smallskip

\noindent Al-Labadi, L. and Evans, M. (2017) Optimal robustness results for
some Bayesian procedures and the relationship to prior-data conflict. Bayesian
Analysis 12, 3, 702-728.\smallskip

\noindent Berger, J.O., Liseo, B. and Wolpert, R.L.$_{.}$ (1999) Integrated
likelihood methods for eliminating nuisance parameters. Statistical Science,
14(1): 1-28. DOI: 10.1214/ss/1009211804\smallskip

\noindent Birnbaum, A. (1964) The anomalous concept of statistical evidence:
axioms, interpretations and elementary exposition. In IMM NYU-332, Courant
Institute of Mathematical Sciences, New York U., New York, USA.\smallskip

\noindent Evans, M. (2015) Measuring Statistical\ Evidence Using Relative
Belief. Chapman and Hall.\smallskip

\noindent Evans, M. and Guo, Y. (2021) Measuring and controlling bias for some
Bayesian inferences and the relation to frequentist criteria. Entropy, 23(2),
190\newline DOI:10.3390/e23020190.\smallskip

\noindent Evans, M. and Jang, G-H. (2011a). Weak informativity and the
information in one prior relative to another. Statistical Science, Vol. 26,
No. 3, 423-439.\smallskip\ 

\noindent Evans, M. and Jang, G-H. (2011b) A limit result for the prior
predictive applied to checking for prior-data conflict. Statistics and
Probability Letters, 81, 1034-1038.\smallskip

\noindent Evans, M. and Moshonov, H. (2006) Checking for prior-data conflict.
Bayesian Analysis, Volume 1, Number 4, 893-914. \smallskip

\noindent Kass, R. E. and Raftery, A. E. (1995). Bayes factors. Journal of the
American Statistical Association, 90, 773-795.\smallskip

\noindent Mulder, J., Berger, J. O.,Pe\~{n}a, V. and Bayarri, M. J. (2017) On
the ubiquity of information inconsistency for conjugate priors.
arXiv:1710.09700\smallskip

\noindent Nott,D., Wang, X., Evans, M., and Englert, B-G. (2020) Checking for
prior-data conflict using prior to posterior divergences. Statistical Science,
35, 2, 234-253.\smallskip

\noindent Popper, K. (1968) The Logic of Scientific Discovery. Harper
Torchbooks.\smallskip

\noindent Royall, R. (1997) Statistical Evidence: A Likelihood Paradigm. CRC
Press.\smallskip

\noindent Rudin, W. (1974) Real and Complex Analysis, Second Edition.
McGraw-Hill.\smallskip

\noindent Salmon, W. (1973) Confirmation. Scientific American, 228,
75--81.\smallskip

\noindent Sarafoglou, A., Barto\v{s}, F., Stefan, A.M., Haaf, J.M., and
Wagenmakers, E.-J. (2022) \textquotedblleft This behavior strikes us as
ideal\textquotedblright: assessment and anticipations of Huisman (2022). DOI:
10.31234/osf.io/rhtcp\smallskip

\noindent Stanford Encyclopedia of Philosophy (2020). Confirmation. In

https://plato.stanford.edu/.
\end{document}